\documentclass[11pt]{amsart}

\usepackage{amsmath,amssymb,amscd}
\usepackage[dvips]{graphics}

\newtheorem{theorem}{Theorem}[section]
\newtheorem{corollary}[theorem]{Corollary}

\newtheorem{lemma}[theorem]{Lemma}
\newtheorem{proposition}[theorem]{Proposition}
\newtheorem{remark}[theorem]{Remark}
\newtheorem{example}[theorem]{Example}
\newtheorem{definition}[theorem]{Definition}

\def\RR{\mathbb{R}}
\def\CC{\mathbb{C}}

\def\PP{\mathbb{P}}
\def\ZZ{\mathbb{Z}}

\def\GL{\operatorname{GL}}
\def\M{\operatorname{Mat}}
\def\adj{\operatorname{adj}}
\def\Id{\operatorname{Id}}
\newcommand{\Span}[1]{\left< #1 \right>}

\begin{document}

\title{Determinantal representations of smooth cubic surfaces}
\author{Anita Buckley and Toma\v z Ko\v sir}

\address{Department of Mathematics, University of Ljubljana,
Jadranska 19, 1000 Ljubljana, Slovenia, 
{\rm e-mail:} {\it anita.buckley@fmf.uni-lj.si,  tomaz.kosir@fmf.uni-lj.si}}%

\begin{abstract}
For every smooth (irreducible) cubic surface $S$ we give an explicit construction of a representative for each of the $72$ 
equivalence classes of determinantal representations. Equivalence classes (under $\GL_3\times \GL_3$ action by left and right 
multiplication) of determinantal representations are in one to one 
correspondence with the sets of six mutually skew lines on $S$ and with the $72$ (two-dimensional) linear systems of twisted 
cubic curves on $S$. Moreover, if a determinantal representation $M$ 
corresponds to lines $(a_1,\ldots,a_6)$ then its transpose $M^t$ corresponds to lines $(b_1,\ldots,b_6)$ which together 
form a Schl\"{a}fli's double-six $a_1\ldots a_6 \choose b_1\ldots b_6$. We also discuss the existence of self-adjoint and definite
determinantal representation for smooth real cubic surfaces. The number of these representations depends on the
Segre type $F_i$. We show that a surface of type $F_i$, $i=1,2,3,4$ has exactly $2(i-1)$ nonequivalent self-adjoint 
determinantal representations none of which is definite, while a surface of type $F_5$ has $24$ nonequivalent 
self-adjoint determinantal representations, $16$ of which are definite.
\end{abstract}

\maketitle

\section{Introduction}

In this paper we study determinantal representations of nonsingular (smooth and irreducible) cubic surfaces. Our first 
objective is to explicitly construct the correspondence in the following classical result which is, for instance, stated 
in \cite[Cor. 6.4]{beauville}. This construction originates in Clebsch \cite{clebsch}. It was probably known already in 
the 19th century, but it is difficult to find a modern reference.
\begin{theorem}\label{thm1}
A smooth cubic surface $S$ in $\PP^3$ allows exactly $72$ nonequivalent determinantal representations. 
There is a one-to-one correspondence between:
\begin{itemize}
\item equivalence classes of determinantal representations of $S$,
\item linear systems of twisted cubic curves on $S$,
\item sets of six lines on $S$ that do not intersect each other.
\end{itemize}
\end{theorem}

Our second objective is to use the explicit construction and study the existence and number of self-adjoint
determinantal representations of smooth real cubic surfaces. Such a representation exists if and only if there is a mutually
self-conjugate double-six on $S$. Our main results are the following: Depending on the Segre type 
$F_i$, $i=1,2,3,4,5$, a smooth real cubic surface has $0,2,4,6$ or $24$, respectively, nonequivalent self-adjoint 
determinantal representations. Surfaces of types $F_1,\ldots,F_4$ have no definite determinantal representations, while a 
surface of type $F_5$ has $16$ definite determinantal representations.

The theory of cubic surfaces has a long history. (See \cite{dolgachev} for a brief outline of the history.) 
It is known since 1849 \cite{cayley} that a smooth cubic surface 
contains $27$ lines. This discovery is one of the first results on surfaces of higher degree and is considered by 
many as the start of modern algebraic geometry. Many mathematicians contributed to the understanding of rich geometry
of smooth cubic surfaces. We refer to Henderson \cite{henderson} or Reid \cite{miles} for the geometry of the $27$ lines
and to Hilbert and Cohn-Vossen  \cite[\S 25]{hilbert} for the geometry of a Schl\"{a}fli's double-six.
It may be hard to find modern references on the geometry of cubic surfaces (an exception being Dolgachev's lecture notes 
\cite{dolgachevwww} available on the web) and many of the results have been known in the 19th century 
\cite{clebsch, cre1, scroter}. A great source for the geometry of real cubic surfaces is Segre \cite{segre}. 

The history of determinantal representations is shorter. Still, Dixon \cite{dixon} showed in 1902 that every curve has a 
determinantal representation. Dickson \cite{dickson} considered determinantal representations of general varieties 
and showed, among other, that every smooth cubic surface has a determinantal representation. Later, determinantal 
representations in general where studied by many authors, for
instance \cite{beauville,cook,room}, and those of cubic surfaces for instance by 
\cite{logar,logar2,geramita}. 
Recently, Vinnikov gave a complete description of equivalence classes of determinantal representations 
for smooth curves \cite{vinnikov2} and of self-adjoint and definite determinantal representations of smooth real curves
\cite{vinnikov,vinnikov3}. Brundu and Logar \cite{logar,logar2} considered singular cubic surfaces 
and showed that all, except those containing one line only, have determinantal representations.

Our main motivation to study determinantal representations comes from possible application to multiparameter
spectral theory \cite{KosirMST}. Classification of (selfadjoint) determinantal representations of algebraic 
hypersurfaces is equivalent to simultaneous classification of tuples of (selfadjoint) matrices. Classically, 
determinantal representations were used to study the rationality of varieties. Vinnikov \cite{vinnikov2, vinnikov3}
stressed another very important motivation to study self-adjoint and definite determinantal representations.
They appear as determinantal representations of discriminant varieties in the theory of commuting nonselfadjoint 
operators in a Hilbert space \cite{LKMV}.

To conclude the introduction we give a brief outline of the paper. In the second section we review the geometry of the 
$27$ lines that we use in the construction. In the third section we give the explicit construction for the 
correspondence between a set of six skew lines on a cubic surface and a determinantal representation in a special form, 
which we call a representation of type $\Re$. In the fourth section we construct the correspondence between equivalence classes of
determinantal representations and systems of twisted cubic curves on $S$. As an example we consider the Fermat cubic surface. 
In sections $5$ and $6$ we present our results on self-adjoint and definite 
determinantal representations of smooth real cubic surfaces.

\section{Classical projective geometry of the 27 lines}

Let $k$ be an algebraically closed field and $S$ a projective cubic surface defined by a smooth
polynomial $F(z_0,z_1,z_2,z_3)$ of degree $3$ over $k.$ 
We will briefly state some classical results in order to introduce the notation.

The most elegant way to study curves on $S$ (our particular interest will be in lines and twisted cubics) is by 
defining $S$ as a blow up of 6 points in the plane, no three collinear and not on a conic. Every nonsingular cubic 
surface in $\PP^3$ can be obtained 
this way. Namely, the linear system of plane cubic curves through assigned points $P_1,\ldots,P_6$ defines an embedding of 
$S$ in $\PP^3.$ 
This way the study of curves on $S$ reduces to study of certain plane curves. The $27$ lines on $S$ are then
the following (see e.g. \cite[Theorem V.4.8.]{hartshorne}): \label{linessev}
\begin{itemize}
\item $a_1,\ldots,a_6$ are the exceptional lines, 
\item $c_{ij}$ is the strict transform of the line through $P_i$ and $P_j$ in $\PP^2,$ where $1\leq i<j\leq 6.$
\item $b_1,\ldots,b_6,$ with $b_j$ being the strict transform of the plane conic through the five $P_i,\ i\neq j.$
\end{itemize}
Observe that $a_1,\ldots,a_6$ are mutually skew, $b_1,\ldots,b_6,$ are mutually skew and $a_i$ intersects $b_j$ 
if and only if $i\neq j.$ Every configuration of $12$ lines on $S$ with this property is called a \textit{Schl\"{a}fli's double-six.}
Every smooth cubic surface $S$ contains $36$ double-sixes of lines.
The $27$ lines have a high degree of symmetry: for 
any set $l_1,\ldots,l_6$ of mutually skew lines on $S$ there exist 
$6$ points in $\PP^2$ and a blow-up for which $l_1,\ldots,l_6$ are the exceptional lines. 
Proof of this can be found in~\cite[Proposition V.4.10.]{hartshorne}. These lines then 
uniquely determine another set of $6$ mutually skew lines to form together a double-six.
Using the above notation the double-sixes on $S$ are:
$$
\left(\begin{array}{ccc}
a_1 & \ldots & a_6 \\
b_1 & \ldots & b_6
\end{array}\right),\ 
\left(\begin{array}{cccccc}
a_i & b_i & c_{kl} & c_{km} & c_{kn} & c_{kp} \\
a_k & b_k & c_{il} & c_{im} & c_{in} & c_{ip}
\end{array}\right),\ 
\left(\begin{array}{cccccc}
a_i & a_k & a_{l} & c_{mn} & c_{mp} & c_{np} \\
c_{kl} & c_{il} & c_{ik} & b_{p} & b_{n} & b_{m}
\end{array}\right).$$
Here $i,k,l,m,n,p$ are all distinct. Let us also mention that a {\em half of a double six}, i.e., a set of six lines from any three
columns of a double-six, belongs to a uniquely determined double-six. This can be easily seen by checking the list of 
all double-sixes.
 
There are many more interesting properties of the configuration of $27$ lines. A plane 
intersecting $S$ in three lines is called a \textit{tritangent plane.} The plane \label{pagepi} 
\begin{equation}
\pi_{ij}=\Span{b_i,a_j}\label{piij}
\end{equation} 
spanned by concurrent lines $b_i,a_j$ intersects the cubic surface $S$ in another line $c_{ij}.$ Note that 
$\pi_{ji}=\Span{b_j,a_i}$ intersects $S$ in the same line $c_{ij}.$ (See for instance \cite[\S 7]{miles}.) 
Every line on $S$ lies exactly on 5 tritangent planes.
In the sequel we 
will use the same symbol for the plane and for its defining linear form. So, the plane $\pi$ is defined by the 
equation $\pi=0$.  

\section{Determinantal representations versus sixes of skew lines}

A \textit{determinantal representation} of $S$ is a $3\times 3$ matrix of linear forms
$$M=M(z_0,z_1,z_2,z_3)=z_0 M_0+z_1 M_1+z_2 M_2+z_3 M_3$$ 
satisfying 
$${\det}M = c F,$$
where $M_0,M_1,M_2,M_3\in \M_3(k)$ and $c \in k,\ c\neq 0.$

Every nonsingular cubic surface $S$ has a determinantal representation. 
We will use the relation between a determinantal representation and six points of the blow-up that
gives $S$, which can be found in the proof of the {existence} of a determinantal representation $M$ in~\cite{geramita}.
Consider the $4$-dimensional linear system of plane cubic curves through $6$ assigned points 
$P_i=(\zeta_i,\eta_i,\xi_i),\ i=1,\ldots ,6.$ A basis 
$F_1,F_2,F_3,F_4$ of the linear system defines a blow-up of $\PP^2(x_0,x_1,x_2).$ 
By the Hilbert-Burch theorem there exists a $3\times 4$ matrix $L$ of linear forms in $x_0,x_1,x_2$ 
whose minors are $F_1,F_2,F_3,F_4$. Note that $L$ has generically rank $3$ and it has rank $2$ exactly 
at points $P_1,\ldots,P_6.$
Define a $3\times 3$ linear matrix $M$ in variables $z_0,z_1,z_2,z_3$ by 
\begin{eqnarray} \label{enacba1}
M\cdot \left(\begin{array}{c}
x_0  \\
x_1 \\
x_2 
\end{array}\right)& = & L\cdot\left(\begin{array}{c}
z_0  \\
z_1 \\
z_2 \\
z_3
\end{array}\right).\end{eqnarray}
One can easily check that $M$ is a determinantal representation of $S.$ 

Conversely, to a given determinantal representation $M$ we can assign $L$ by (\ref{enacba1}) above.
The rank of $L$ in 
$P=(\zeta,\eta,\xi)\in \PP^2$ equals $2$ if and only if the three planes in $\PP^3$ with equations
$$M \cdot \left(\begin{array}{c}
\zeta\\
\eta\\
\xi
\end{array}\right)=\left(\begin{array}{c}
0\\
0\\
0
\end{array}\right)$$
intersect in a line. Note that the lines obtained this way are exactly the exceptional lines of the 
blow-up~\cite{geramita}. They are mutually skew and we call them the \textit{six skew lines corresponding 
to determinantal representation} $M.$

Two determinantal representations $M$ and $M'$ are \textit{equivalent} if there exist $X,Y\in \GL_3(k)$
such that
$$M'=XMY.$$

\begin{lemma}\label{conv} 
If the determinantal representations $M$ and $M'$ are equivalent then the corresponding sets of $6$ skew lines are equal.
\end{lemma}

\begin{proof} Let $a_1,\ldots,a_6$ be the lines corresponding to $M$ over points 
$P_i=(\zeta_i,\eta_i,\xi_i),\ i=1,\ldots,6$ 
as above and let $M'=XMY$ with $X,Y\in\GL_3(k).$ Then 
$$\mathcal{L}in \left\{ M'Y^{-1}\left(\begin{array}{c}
\zeta_i\\
\eta_i \\
\xi_i
\end{array}\right)\right\}=\mathcal{L}in \left\{XM\left(\begin{array}{c}
\zeta_i\\
\eta_i \\
\xi_i
\end{array}\right)\right\}=\mathcal{L}in \left\{M\left(\begin{array}{c}
\zeta_i\\
\eta_i \\
\xi_i
\end{array}\right)\right\}$$ 
and so the line over $P_i$ for $M$ is the same as the line over $Y^{-1}P_i$ for $M'$. 
\end{proof}

In the rest of the section we will explicitly construct
all (up to equivalence) determinantal representations from the $27$ lines on $S.$ 

We fix an arbitrary set of six mutually skew lines on $S.$ 
Without loss of generality name them $a_1,\ldots,a_6$. There is a unique set of $6$
skew lines on $S$ such that 
$$\left(\begin{array}{ccc}
a_1 & \ldots & a_6 \\
b_1 & \ldots & b_6
\end{array}\right)$$
is a double-six. 
Consider the tritangent planes 
$$\pi_{12},\pi_{23},\pi_{31},\pi_{13},\pi_{21},\pi_{32} $$
defined in (\ref{piij}). 
On $S$ they intersect in $9$ lines 
$$a_1,a_2,a_3,b_1,b_2,b_3,c_{12},c_{13},c_{23}$$
which together with another point on $S$ outside these lines determine $F.$ 
This follows from the fact that cubic surfaces are parametrized by the points in $\PP^{19}$,
whereas $9$ lines and a point induce $2\cdot 9+1=19$ conditions. Then $F=\pi_{12}\pi_{23}\pi_{31}+
\lambda\pi_{13}\pi_{21}\pi_{32}$ for some nonzero $\lambda \in k$, and we may without loss assume that $\lambda=1$.
Thus 
\begin{equation}\label{R}
F=\pi_{12}\pi_{23}\pi_{31}+\pi_{13}\pi_{21}\pi_{32}\ \mbox{ and }\ \Re=\left(\begin{array}{ccc}
0 & \pi_{12} & \pi_{13}  \\
\pi_{21} & 0 & \pi_{23}\\
\pi_{31} & \pi_{32} & 0
\end{array}\right).
\end{equation}
Note that $\Re$ is a determinantal representation of $S$. We say that a determinantal representation is {\em of type} 
$\Re$ if it is in the above form. Similar form was considered in \cite[\S 1.6]{kollar}.

\begin{proposition} \label{doublrem}
The two sets of six skew lines corresponding to $\Re$ and to its transpose $\Re^t$ form the double-six
$a_1\ldots a_6 \choose b_1\ldots b_6$, which was our construction's input data.
\end{proposition}

\begin{proof}
From the construction it is obvious that the intersections of planes corresponding to
$$\Re \left(\begin{array}{c}
1\\
0 \\
0
\end{array}\right),\ \Re \left(\begin{array}{c}
0\\
1 \\
0
\end{array}\right),\ \Re \left(\begin{array}{c}
0\\
0 \\
1
\end{array}\right)$$
are $a_1,a_2,a_3$, respectively, and those corresponding to
$$\Re^t \left(\begin{array}{c}
1\\
0 \\
0
\end{array}\right),\  \Re^t \left(\begin{array}{c}
0\\
1 \\
0
\end{array}\right),\ \Re^t \left(\begin{array}{c}
0\\
0 \\
1
\end{array}\right)$$
are $b_1,b_2,b_3$, respectively. 

There exist three more points $(\zeta_k,
\eta_k,
\xi_k)\in \PP^2,\ k=4,5,6$
such that the intersections of planes given by 
$$\Re  \left(\begin{array}{c}
\zeta_4\\
\eta_4 \\
\xi_4
\end{array}\right), \Re  \left(\begin{array}{c}
\zeta_5\\
\eta_5 \\
\xi_5
\end{array}\right),\Re  \left(\begin{array}{c}
\zeta_6\\
\eta_6\\
\xi_6
\end{array}\right)$$
induce skew lines $l_4,l_5,l_6$, respectively.

Observe that 
$$a_2 \subset \pi_{12},\ a_3 \subset \pi_{13},\
l_4 \subset \eta_4 \pi_{12}+\xi_4 \pi_{13},\ l_5 \subset \eta_5 \pi_{12}+\xi_5 \pi_{13},\
l_6 \subset \eta_6 \pi_{12}+\xi_6 \pi_{13}$$ 
and that $b_1$ lies on all these tritangent planes. These planes are all different since they contain different lines
that are skew. 
Because every two lines in a plane intersect, $b_1$ intersects $a_2, a_3,l_4,l_5,l_6$ and is disjoint with $a_1.$
In the same way we prove that 
$b_2$ intersects $a_1, a_3,l_4,l_5,l_6$, but not $a_2$, and that
$b_3$ intersects $a_1, a_2,l_4,l_5,l_6$, but not $a_3$.  

Applying the same arguments to $\Re^t$ yields $6$ mutually skew lines $b_1,b_2,b_3,l'_4,l'_5,l'_6$
and relations:
$a_1$ intersects $b_2, b_3,l'_4,l'_5,l'_6,$, but not $b_1$,\ $a_2$ intersects $b_1, b_3,l'_4,l'_5,l'_6,$ but not 
$b_2$, \ $a_3$ intersects $b_1, b_2,l'_4,l'_5,l'_6,$ but not $b_3$. 
These relations are fulfilled if and only if the above lines form a double-six and 
the remaining lines $l_4,l_5,l_6$ and $l_4',l'_5,l'_6$ have to be equal to $a_4,a_5,a_6$ and $b_4,b_5,b_6$, respectively. 
(See Section 2.)
\end{proof} 

We have set everything needed to prove the converse of Lemma~\ref{conv}.
We will show that every determinantal representation of $S$ is equivalent to one of type $\Re$. 

\begin{proposition} \label{mainth} If determinantal representations $M$ and $M'$ of $S$ define the same set of six
skew lines then they are equivalent. 
\end{proposition}

\begin{proof} Let $a_1,\ldots,a_6$ be the set of lines corresponding to $M$ over 
the points $P_i=(\zeta_i,\eta_i,\xi_i),\ i=1,\ldots ,6.$. It suffices to prove that $M$ is equivalent
to the representation $\Re$. We will prove the statement by explicitly constructing the equivalence.  

Consider representation 
$$M'=M\left(\begin{array}{ccc}
\zeta_1 & \zeta_2 & \zeta_3 \\
\eta_1 &\eta_2 & \eta_3 \\
\xi_1 & \xi_2 & \xi_3
\end{array}\right),$$
which is equivalent to $M$.
Since the three planes in $\PP^3$ given by linear equations 
$$M' \left(\begin{array}{c}
1\\
0 \\
0
\end{array}\right)=\left(\begin{array}{c}
m'_{11}\\
m'_{21}\\
m'_{31}
\end{array}\right)=\left(\begin{array}{c}
0\\
0\\
0
\end{array}\right)$$
intersect in a line ($a_1$ actually), there exists a nontrivial linear combination 
$$\alpha_1 m'_{11}+\beta_1 m'_{21} +\gamma_1 m'_{31}=0,$$
i.e., $\alpha_1,\beta_1,\gamma_1\in k$ are not all equal to zero.
After multiplying $M'$ from the left by an invertible constant matrix 
$$A_1=\left(\begin{array}{ccc}
\alpha_1 & \beta_1 & \gamma_1 \\
\ast &\ast & \ast \\
\ast & \ast & \ast
\end{array}\right)$$
we obtain 0 in the ${(1,1)}$th entry. 

Observe that the ${(1,2)}$th and ${(1,3)}$th entries of $A_1 M'$ are nonzero since $F$ is irreducible. 
The planes in the second column of $A_1 M'$ again intersect in a line, namely $a_2.$ 
Therefore there exists a nontrivial combination of these planes with coefficients
$(\alpha_2, \beta_2, \gamma_2)$ different from a multiple of $(1,0,0)$.
By multiplying $A_1 M'$ from the left by an invertible matrix
$$A_2=\left(\begin{array}{ccc}
1 & 0&0 \\
\alpha_2 & \beta_2 & \gamma_2 \\
\ast & \ast & \ast
\end{array}\right)$$
we obtain $0$ on the ${(2,2)}$th place.

Observe that the ${(1,3)}$th and ${(2,3)}$th entry of $A_2 A_1 M'$ are not proportional, and thus both nonzero. 
By the same argument as above, there exists a linear combination with coefficients
$(\alpha_3, \beta_3, \gamma_3)$
of the entries of third column of $A_2 A_1 M'$ with $\gamma_3\neq 0.$

Finally, 
$$A_3=\left(\begin{array}{ccc}
1 & 0&0 \\
0&1&0 \\
\alpha_3 & \beta_3 & \gamma_3
\end{array}\right)$$
is such that 
$$M''=A_3 A_2 A_1 M'=\left(\begin{array}{ccc}
0 & m''_{12}& m''_{13} \\
m''_{21}&0&m''_{23} \\
m''_{31} & m''_{32} &0
\end{array}\right)$$
has zeros along the diagonal.

Since $M''$ is equivalent to $M'$ the intersections of planes given by the entries of
$$M''  \left(\begin{array}{c}
1\\
0 \\
0
\end{array}\right),\ M''  \left(\begin{array}{c}
0\\
1 \\
0
\end{array}\right),\ M''  \left(\begin{array}{c}
0\\
0 \\
1
\end{array}\right)$$
induce the lines $a_1,a_2,a_3$, respectively, and there exist three more points in $\PP^2$,
namely
$$\left(\begin{array}{c}
\zeta_k''\\
\eta_k'' \\
\xi_k''
\end{array}\right)=\left(\begin{array}{ccc}
\zeta_1 & \zeta_2 & \zeta_3 \\
\eta_1 &\eta_2 & \eta_3 \\
\xi_1 & \xi_2 & \xi_3
\end{array}\right)^{-1} \left(\begin{array}{c}
\zeta_k\\
\eta_k \\
\xi_k
\end{array}\right),\ k=4,5,6,$$
such that the intersections of planes corresponding to
$$M''  \left(\begin{array}{c}
\zeta_4''\\
\eta_4'' \\
\xi_4''
\end{array}\right),\ M''  \left(\begin{array}{c}
\zeta_5''\\
\eta_5'' \\
\xi_5''
\end{array}\right),\ M''  \left(\begin{array}{c}
\zeta_6''\\
\eta_6'' \\
\xi_6''
\end{array}\right)$$
induce $a_4,a_5,a_6$, respectively.\\

Consider now the transpose representation $(M'')^t$ and its first corresponding line $l$ defined as
intersection of the planes corresponding to the entries of
$$(M'')^t \left(\begin{array}{c}
1\\
0 \\
0
\end{array}\right)=\left(\begin{array}{c}
0\\
m''_{12} \\
m''_{13}
\end{array}\right).$$
Observe that 
$$a_2 \subset m''_{12}, a_3 \subset m''_{13},
a_4 \subset \eta_4'' m''_{12}+\xi_4'' m''_{13},a_5 \subset \eta_5'' m''_{12}+\xi_5'' m''_{13},
a_6 \subset \eta_6'' m''_{12}+\xi_6'' m''_{13}$$ 
and that $l$ lies on all these tritangent planes. These planes are all different since they contain different 
lines that are skew. 
Since every line of $S$ lies on $5$ tritangent planes it follows that $l\cap a_1=\emptyset$, otherwise $l$ and 
$a_1$ would span another tritangent plane.
Because every two lines in a plane intersect, $l$ intersects $a_2, a_3,a_4,a_5,a_6$ and is disjoint with $a_1.$
Therefore $l=b_1.$ 
We continue in the same way to prove that 
$$\mbox{the line defined by }(M'')^t \left(\begin{array}{c}
0\\
1 \\
0
\end{array}\right) \mbox{ is }b_2$$
 and 
$$\mbox{the line defined by }(M'')^t \left(\begin{array}{c}
0\\
0 \\
1
\end{array}\right)\mbox{ is }b_3.$$

We claim that $m_{ij}''$ yields an equation of the tritangent plane $\pi_{ij}.$
Indeed, the plane $m_{ij}''$ contains the lines $b_i$ and $a_j.$ The incidence relation can be read from the diagram
$$\begin{array}{c}
\left(\begin{array}{c}
0   \\
m''_{21}\\
m''_{31}
\end{array}\right. \\
\downarrow \\
a_1 
\end{array}
\begin{array}{c}
\begin{array}{c}
m''_{12}   \\
0\\
m''_{32}
\end{array} \\
\downarrow \\
a_2 
\end{array}
\begin{array}{c}
\left.\begin{array}{c}
m''_{13}  \\
m''_{23}\\
0
\end{array}\right) \\
\downarrow \\
a_3 
\end{array}
\begin{array}{c}
\begin{array}{c}
\rightarrow   \\
\rightarrow  \\
\rightarrow 
\end{array} \\
\\
\\
\end{array}
\begin{array}{c}
\begin{array}{c}
b_1   \\
b_2 \\
b_3
\end{array} \\
\\
\\
\end{array}
$$
This means that there exist nonzero $\alpha,\alpha',\alpha'',\beta,\beta',\beta''$ in $k$ such that
$\alpha\alpha'\alpha''=\beta\beta'\beta''$
and
$$\left(\begin{array}{ccc}
0 & \pi_{12} & \pi_{13}  \\
\pi_{21} & 0 & \pi_{23}\\
\pi_{31} & \pi_{32} & 0
\end{array}\right)=\left(\begin{array}{ccc}
0 & \alpha m''_{12}& \beta m''_{13} \\
\beta' m''_{21} & 0 &\alpha' m''_{23} \\
\alpha'' m''_{31} & \beta'' m''_{32} &0
\end{array}\right),$$
which is equal to
$$
\left(\begin{array}{ccc}
\frac{1}{\alpha'} & 0 & 0 \\
0 &   \frac{1}{\beta} & 0\\
0 & 0 & \frac{\alpha''}{\beta \beta'}
\end{array}\right)
\left(\begin{array}{ccc}
0 & m''_{12}& m''_{13} \\
m''_{21}&0&m''_{23} \\
m''_{31} & m''_{32} &0
\end{array}\right)
\left(\begin{array}{ccc}
\beta \beta' & 0 & 0 \\
0 &  \alpha \alpha' & 0\\
0 & 0 & \alpha' \beta
\end{array}\right).$$
We proved that
$$\Re=\left(\begin{array}{ccc}
\frac{1}{\alpha'} & 0 & 0 \\
0 &   \frac{1}{\beta} & 0\\
0 & 0 & \frac{\alpha''}{\beta \beta'}
\end{array}\right)
A_3 A_2 A_1 M\left(\begin{array}{ccc}
\zeta_1 & \zeta_2 & \zeta_3 \\
\eta_1 &\eta_2 & \eta_3 \\
\xi_1 & \xi_2 & \xi_3
\end{array}\right)
\left(\begin{array}{ccc}
\beta \beta' & 0 & 0 \\
0 &  \alpha \alpha' & 0\\
0 & 0 & \alpha' \beta
\end{array}\right),$$
which shows that $M$ is equivalent to the representation $\Re$ of (\ref{R}).
\end{proof}

\begin{remark}{\rm The proof of Proposition~\ref{mainth}  describes a method how to explicitly construct 
$X,Y\in\GL_3(k)$ such that $XMY$ is equal to~$\Re$.} 

{\rm Note also that Lemma \ref{conv} and Proposition \ref{mainth} prove that there is a bijective correspondence
between the equivalence classes of determinantal representations of $S$ and sets of six skew lines on $S$. The proof
provides an explicit construction of the correspondence. }
\end{remark}

\begin{corollary}\label{longlem}
If the lines $a_1,\ldots,a_6$ correspond to $M$ and the lines $b_1,\ldots,b_6$ correspond to $M^t,$ then 
$a_1,\ldots,a_6$ and $b_1,\ldots,b_6$ form a double-six.
\end{corollary}

\begin{proof}
We proved that $M$ is equivalent to $\Re$ because they induce the same 6 lines $a_1,\ldots,a_6$. Then $M^t$ and
$\Re^t$ are also equivalent and by Lemma~\ref{conv} induce the same set of lines $b_1,\ldots,b_6.$ 
By Proposition~\ref{doublrem} they form together a double-six.
\end{proof}

\begin{corollary}\label{symm}
A smooth cubic surface can not have a symmetric determinantal representation.
\end{corollary}

\section{Determinantal representations and twisted cubic curves}

So far we proved and explicitly constructed correspondence between equivalence classes 
of determinantal representations of $S$ and sixes of skew lines on $S$. More precisely, there is 1-1 
correspondence between pairs $M,\,M^t$ and double-sixes.

In this section we find 1-1 correspondence between equivalence classes 
of determinantal representations and linear systems of twisted cubics on $S$. \\

The notion of equivalent representations is strongly connected to the equivalence of the following line bundles. 
If $M$ is a determinantal representation for $S,$ define the \textit{corresponding vector bundle} $\epsilon$ 
along $S$ by
$$\epsilon={\ker}M.$$
Since $S$ is nonsingular, $\epsilon$ is a line bundle by
\begin{lemma} \label{dimone}
The corresponding vector bundle of a determinantal representation of a smooth surface is of rank one. 
\end{lemma}
\begin{proof} 
Let 
$$M(z_0,z_1,z_2,z_3)=\left(m_{ij}^0 z_0 +m_{ij}^1 z_1 +m_{ij}^2 z_2 +m_{ij}^3 z_3\right)_{i,j}$$ 
satisfy $\mbox{det}M = c F$
and denote by $A_{ij}$ the $(i,j)$th $2$nd-order cofactor of $M.$ 
Using the formula for the differentiation of a determinant 
and row expansion we get
$$\frac{\partial F}{\partial z_k}=\frac{1}{c}\sum_{i,j}m_{ij}^k A_{ij.}$$
At every regular point of $F$ not all partial derivatives are 0 and therefore not all $2\times 2$ minors
of $M$ are 0. Thus $M$ has rank 2.
\end{proof}

Next we introduce the adjoint matrix $$\widetilde{M}={\adj}M,$$ whose $(i,j)$th entry is the $(j,i)$th cofactor 
of $M$. From the identity
$$M \widetilde{M}=({\det}M)\, \Id=(c\, F)\,\Id $$ 
and Lemma~\ref{dimone} we immediately deduce that
\begin{itemize}
\item every entry of $\widetilde{M}$ is a homogeneous polynomial of degree 2,
\item matrix $\widetilde{M}$ has rank 1 at every point on $S.$
\end{itemize}
Thus on $S$ every column of $\widetilde{M}$ is in ${\ker}M$ and its zero locus, i.e. the zero locus defined 
by the entries of a column, is a twisted cubic curve on $S$. (See e.g. \cite[Lecture 12]{harris}.)
Since $\det M\neq 0$ we get three different twisted cubic curves corresponding to the columns of $\widetilde{M}$ 
(which are equivalent as effective divisors on $S$). 
In other words, the columns define a basis of a 2-dimensional linear system of twisted cubic curves on $S$. 
We call any of twisted cubic curves obtained above a \textit{twisted cubic corresponding to determinantal 
representation $M$.}

Recall the identification~\cite{shaf} of line bundles with invertible sheaves. 
In other words, we could define $\epsilon$ corresponding to $M$ by an exact sequence
\begin{eqnarray}\label{enac2}
0 \longrightarrow \mathcal{O}_{\PP^3}(-1)^3 \stackrel{M}{\longrightarrow} \mathcal{O}^3_{\PP^3} 
\longrightarrow \epsilon \longrightarrow 0.
\end{eqnarray}
The invertible sheaf $\epsilon$ defined in such a way is obviously supported on $S$.
By~\cite[Theorem 1.1.]{cook} matrices $M$ and $M'$ are equivalent if and only if the corresponding sheaves 
$\epsilon$ and $\epsilon'$ are isomorphic. 
Note that (\ref{enac2}) induces long exact sequence with 
$$\mbox{H}^0(\PP^3,\mathcal{O}_{\PP^3}(-1)^3  )=\mbox{H}^1(\PP^3,\mathcal{O}_{\PP^3}(-1)^3 )=0$$
and hence 
$$\mbox{H}^0(\PP^3,\mathcal{O}^3_{\PP^3})\cong \mbox{H}^0(S,\epsilon)$$
has dimension 2.

Moreover, Beauville~\cite[Proposition 6.2.]{beauville} proves that 
the cokernel of 
$\mathcal{O}_{\PP^3}(-1)^3 \stackrel{M}{\longrightarrow} \mathcal{O}^3_{\PP^3} $ is isomorphic to
$\mathcal{O}_S(C),$ where $C$ is a smooth projectively normal curve on $S$ of degree 3 and genus 0. 
Thus $C$ is a twisted cubic on $S$ as we have already seen. By~\cite[Proposition II.6.13.]{hartshorne} 
$\mathcal{O}_S(C)\cong \mathcal{O}_S(C')$ 
if and only if $C$ and $C'$ are equivalent as effective divisors on $S.$
We can conclude that our linear system of twisted cubic curves is in fact a complete linear 
system $\mbox{H}^0(S,\mathcal{O}_S(C))$ of a twisted cubic curve $C$ corresponding to $M$.

\begin{lemma}
Equivalent determinantal representations induce the same linear system of twisted cubic curves.
\end{lemma}

\begin{proof}
Let $M$ be a determinantal representation of $S$. The columns of the adjoint matrix $\widetilde{M}$
form a basis of the corresponding linear system of 
twisted cubics on $S$. In other words 
$$\widetilde{M}=\left(\begin{array}{ccc}
q_1 & q_1 & q_1 \\
q_2 & q_2 & q_2 \\
q_3 & q_3 & q_3
\end{array}\right)\left(\begin{array}{ccc}
1 & 0 & 0 \\
0 & r & 0 \\
0 & 0 & t
\end{array}\right)$$
where $q_i$ are quadratic polynomials and $r,t$ rational functions on $S.$ 

It is now obvious that an equivalent representation $XMY,\ X,Y\in\GL_3(k)$ induces the same linear system,
just another basis.
\end{proof}

It remains to show that every linear system of twisted cubics on $S$ indeed comes from some representation.

Consider $S$ as a blow-up of $6$ points of the plane and  denote by
$A_1,\ldots,A_6$ the linear equivalence classes of the exceptional lines. Let $L\in\mbox{Pic}\,S$ be
the pull-back of a line in $\PP^2.$ 

\begin{proposition}\cite[Proposition V. 4.8]{hartshorne} \label{divhart} Let $S, A_i$ and $L$ be as above. Then:
\begin{enumerate}
\item $\mbox{Pic}S\cong \ZZ^7$ is generated by $A_1,\ldots,A_6$ and $L;$
\item an effective divisor $D\sim \alpha L-\sum \beta_i A_i$ on $S$ 
has degree 
$$d=3\alpha-\sum \beta_i$$ as a curve in $\PP^3;$
\item the self-intersection of $D$ is
$$D^2=2(p_a(D)-1)+d=\alpha^2-\sum \beta^2_i,$$
where $p_a$ denotes the arithmetic genus of $D$.
\end{enumerate}
\end{proposition}

The following is a well known result. 

\begin{proposition} 
Every smooth cubic surface $S$ contains $72$ linear systems of twisted cubic curves.
\end{proposition}

\begin{proof} We need to count all nonequivalent effective divisors $D\sim \alpha L-\sum \beta_i A_i$  
on $S$ of degree $d=3$ and genus $p_a=0.$ 
We will show that the only integers satisfying these conditions are
$$\begin{array}{ccc}
\alpha & \{\beta_1,\ldots,\beta_6 \} & \mbox{number of solutions} \\
 & & \\
1 & \{0,0,0,0,0,0\} & 1 \\
2 & \{1,1,1,0,0,0\} & 20 \\
3 & \{2,1,1,1,1,0\} & 30 \\
4 & \{2,2,2,1,1,1\} & 20 \\
5 & \{2,2,2,2,2,2\} & 1 
\end{array}$$
Recall Schwarz's inequality for real sequences
$$|\sum x_i y_i|^2\leq |\sum x_i^2|\cdot|\sum y_i^2|$$
and take $x_i=1,\ y_i=\beta_i,\ i=1,\ldots,6.$
By Proposition~\ref{divhart} we find
$$(3\alpha-3)^2=(\sum \beta_i)^2\leq 6\sum \beta_i^2=6(\alpha^2-1).$$
Solving the quadratic equation for $\alpha$ we get $1\leq \alpha \leq 5.$
One quickly obtains all possible integers $\beta_i$ by trial.
\end{proof}

Recall representation $\Re$ of (\ref{R}) corresponding to the double-six $a_1\ldots a_6 \choose b_1\ldots b_6$.
Its adjoint is equal to
$$\begin{array}{ll}
\widetilde{\Re} & =
\left(\begin{array}{ccc}
-\pi_{23}\pi_{32} & \pi_{13}\pi_{32} & \pi_{12}\pi_{23} \\
\pi_{23}\pi_{31}  & -\pi_{13}\pi_{31} & \pi_{13}\pi_{21}\\
\pi_{21}\pi_{32}  & \pi_{12}\pi_{31} & -\pi_{12}\pi_{21}
\end{array}\right) \\
 & =\left(\begin{array}{ccc}
-\pi_{23}\pi_{32} & -\pi_{23}\pi_{32} & -\pi_{23}\pi_{32} \\
\pi_{23}\pi_{31} & \pi_{23}\pi_{31} & \pi_{23}\pi_{31}\\
\pi_{21}\pi_{32} & \pi_{21}\pi_{32} & \pi_{21}\pi_{32}
\end{array}\right)\left(\begin{array}{ccc}
1 & 0 & 0 \\
0 & -\frac{\pi_{13}}{\pi_{23}} & 0 \\
0 & 0 & -\frac{\pi_{12}}{\pi_{32}}
\end{array}\right)
\end{array}$$

The zero locus of the first column is a \textit{degenerated twisted cubic} consisting of lines
$$\pi_{23}\cap \pi_{21}=b_2,\ \pi_{23}\cap \pi_{32}=c_{23},\ \pi_{32}\cap \pi_{31}=b_3.$$

As divisors, the 27 lines on $S$ are: 
$$a_i\sim A_i,\ c_{ij}\sim L-A_i-A_j,\ b_j\sim 2L-\sum_{i\neq j} A_i,$$
This follows from the discussion in \S 2 of $S$ as a blow-up of six points.
Thus the degenerated twisted cubic curve corresponding to the first column of $\Re$ is an effective divisor
$$ 5L-2\sum_{i=1}^6 A_i\sim b_2+b_3+c_{23}.$$
It now follows easily that twisted cubics corresponding to equivalent representations (or the same 
set of six skew lines) are linearly equivalent as divisors on $S$. 
And a different set of six skew lines induces a linearly nonequivalent divisor, i.e. different
linear system of twisted cubic curves.

This proves that on $S$ there are exactly $72$ linear systems of twisted cubics, each corresponding to a set of six
skew lines and concludes a proof of Theorem \ref{thm1}.

\begin{remark} {\rm The two-dimensional linear system $|5L-2\sum_{i=1}^6 A_i|$ corresponding to the six lines 
$a_1,\ldots,a_6,$ contains ${6 \choose 2}=15$ degenerated twisted cubics consisting of 3 lines $b_i,b_j,c_{ij}$.}
\end{remark}


\begin{remark}{\rm
In 1856 Steiner~\cite{henderson},~\cite{dolgachev} introduced sets of 9 lines on $S$ that can be put in a 
$3\times 3$ array such that lines in each row and column generate a tritangent plane. Here is the list of all possible
Steiner sets:
$$\begin{array}{ccc}
a_i & b_j & c_{ij} \\
b_k & c_{jk} & a_j \\
c_{ik} & a_k & b_i
\end{array}\ \ \ \ \begin{array}{ccc}
c_{ij} & c_{kl} & c_{mn} \\
c_{ln} & c_{im} & c_{jk} \\
c_{km} & c_{jn} & c_{il}
\end{array}\ \ \ \ \begin{array}{ccc}
a_{i} & b_{j} & c_{ij} \\
b_{k} & a_{l} & c_{kl} \\
c_{ik} & c_{jl} & c_{mn}
\end{array}
$$
Denote by $\rho_1,\rho_2,\rho_3$ and $\sigma_1,\sigma_2,\sigma_3$ the tritangent planes corresponding 
to rows and columns, respectively, of a Steiner set. These planes 
form a so called \textit{Triederpaar}. There are $120$ Triederpaars and thus $120$ essentially different 
representations of a cubic surface by an equation of the form
$$\rho_1\rho_2\rho_3+\sigma_1\sigma_2\sigma_3=0.$$

\noindent
\textbf{Question:} How do we relate the above $120$ representations to the $72$ classes of nonequivalent 
determinantal representations?\\

In fact, every Triederpaar induces $2\cdot 3!\cdot 3!$
determinantal representations
$$\mathcal{M}=\left(\begin{array}{ccc}
0 & \rho_i & \sigma_p \\
\sigma_n & 0 & \rho_j \\
\rho_k & \sigma_m & 0
\end{array}\right)\ \mbox{ and }\ \mathcal{M}^t$$
for $\{i,j,k\}=\{p,m,n\}=\{1,2,3\}$, which can be obviously partitioned in sets of $6$ equivalent ones 
(obtained by permutation group $S_3$ acting on the same rows and columns). 

Recall that the lines defined by the columns  
$$\mathcal{M}\left(\begin{array}{c}
1\\
0\\
0
\end{array}\right),\ 
\mathcal{M}\left(\begin{array}{c}
0\\
1\\
0
\end{array}\right),\ 
\mathcal{M}\left(\begin{array}{c}
0\\
0\\
1
\end{array}\right),\ 
\mathcal{M}^t\left(\begin{array}{c}
1\\
0\\
0
\end{array}\right),\
\mathcal{M}^t\left(\begin{array}{c}
0\\
1\\
0
\end{array}\right),\
\mathcal{M}^t\left(\begin{array}{c}
0\\
0\\
1
\end{array}\right)$$
form a half {of a double-six} that can be uniquely extended to the whole double-six 
corresponding to the pair of determinantal representations $\mathcal{M}$ and $\mathcal{M}^t.$ 
Determinantal representations with zeros along the diagonal are equivalent if and only if the corresponding
halves of a double-six extend to the same double-six. There are ${6 \choose 3}=20$ ways of choosing
a half of a given double-six. Thus we find that the $120$ Steiner sets give
$$120\cdot\frac{2\cdot 3!\,3!}{6}\cdot\frac{1}{20}=72$$
nonequivalent determinantal representations on $S$.}
\end{remark}

\begin{example}\label{exferm}{\rm Consider Fermat surface $S$ given by the equation
$$F=z_0^3+z_1^3+z_2^3+z_3^3=0.$$ 
We represent a line 
$$\alpha_0 z_0+\alpha_1 z_1+\alpha_2 z_2+\alpha_3 z_3=0,\ \ \beta_0 z_0+\beta_1 z_1+\beta_2 z_2+\beta_3 z_3=0$$
by a $2\times 4$ matrix
$$\left(\begin{array}{cccc}
\alpha_0 & \alpha_1 & \alpha_2 & \alpha_3 \\
\beta_0 & \beta_1 & \beta_2 & \beta_3
\end{array}\right).$$
The 27 lines on $S$ are then: \label{linets}
$$\begin{array}{c}
\left(\begin{array}{cccc}
1 & 1 & 0 & 0\\
0 & 0 & 1 & 1
\end{array}\right),\  
\left(\begin{array}{cccc}
1 & 1 & 0 & 0\\
0 & 0 & 1 & \omega
\end{array}\right),\  
\left(\begin{array}{cccc}
1 & 1 & 0 & 0\\
0 & 0 & \omega & 1
\end{array}\right), \\
\left(\begin{array}{cccc}
1 & \omega & 0 & 0\\
0 & 0 & 1 & 1
\end{array}\right),\
\left(\begin{array}{cccc}
1 & \omega & 0 & 0\\
0 & 0 & 1 & \omega
\end{array}\right),\
\left(\begin{array}{cccc}
1 & \omega & 0 & 0\\
0 & 0 & \omega & 1
\end{array}\right),\\
\left(\begin{array}{cccc}
\omega & 1 & 0 & 0\\
0 & 0 & 1 & 1
\end{array}\right),\  
\left(\begin{array}{cccc}
\omega & 1 & 0 & 0\\
0 & 0 & 1 & \omega
\end{array}\right),\  
\left(\begin{array}{cccc}
\omega & 1 & 0 & 0\\
0 & 0 & \omega & 1
\end{array}\right), \\
\left(\begin{array}{cccc}
1 & 0 & 1 & 0\\
0 & 1 & 0 & 1
\end{array}\right),\  
\left(\begin{array}{cccc}
1 & 0 & 1 & 0\\
0 & 1 & 0 & \omega
\end{array}\right),\  
\left(\begin{array}{cccc}
1 & 0 & 1 & 0\\
0 & \omega & 0 & 1
\end{array}\right), \\
\left(\begin{array}{cccc}
1 & 0 & \omega & 0\\
0 & 1 & 0 & 1
\end{array}\right),\  
\left(\begin{array}{cccc}
1 & 0 & \omega & 0\\
0 & 1 & 0 & \omega
\end{array}\right),\  
\left(\begin{array}{cccc}
1 & 0 & \omega & 0\\
0 & \omega & 0 & 1
\end{array}\right), \\
\left(\begin{array}{cccc}
\omega & 0 & 1 & 0\\
0 & 1 & 0 & 1
\end{array}\right),\  
\left(\begin{array}{cccc}
\omega & 0 & 1 & 0\\
0 & 1 & 0 & \omega
\end{array}\right),\  
\left(\begin{array}{cccc}
\omega & 0 & 1 & 0\\
0 & \omega & 0 & 1
\end{array}\right), \\
\left(\begin{array}{cccc}
1 & 0 & 0 & 1\\
0 & 1 & 1 & 0
\end{array}\right),\  
\left(\begin{array}{cccc}
1 & 0 & 0 & 1\\
0 & 1 & \omega & 0
\end{array}\right),\  
\left(\begin{array}{cccc}
1 & 0 & 0 & 1\\
0 & \omega & 1 & 0
\end{array}\right), \\
\left(\begin{array}{cccc}
1 & 0 & 0 & \omega\\
0 & 1 & 1 & 0
\end{array}\right),\  
\left(\begin{array}{cccc}
1 & 0 & 0 & \omega\\
0 & 1 & \omega & 0
\end{array}\right),\  
\left(\begin{array}{cccc}
1 & 0 & 0 & \omega\\
0 & \omega & 1 & 0
\end{array}\right), \\
\left(\begin{array}{cccc}
\omega & 0 & 0 & 1\\
0 & 1 & 1 & 0
\end{array}\right),\  
\left(\begin{array}{cccc}
\omega & 0 & 0 & 1\\
0 & 1 & \omega & 0
\end{array}\right),\  
\left(\begin{array}{cccc}
\omega & 0 & 0 & 1\\
0 & \omega & 1 & 0
\end{array}\right),
\end{array}
$$
where $\omega$ is a primitive third root of unity.
Consider determinantal representation
\begin{equation}\label{Fermat1}
M=\left(\begin{array}{ccc}
0 & z_0+z_1 & z_2+z_3 \\
\omega z_2+z_3 & 0 & z_0+\omega z_1 \\
\omega z_0+z_1 & z_2+\omega z_3 & 0
\end{array}\right).
\end{equation} 
Recall from equality~(\ref{enacba1}) that
\begin{eqnarray} \label{enacba7}
M\cdot \left(\begin{array}{c}
x_0  \\
x_1 \\
x_2 
\end{array}\right)& = & 
\left(\begin{array}{cccc}
x_1 & x_1 & x_2 & x_2\\
x_2 & \omega x_2 & \omega x_0 & x_0 \\
\omega x_0 & x_0 & x_1 & \omega x_1
\end{array}\right)
\cdot\left(\begin{array}{c}
z_0  \\
z_1 \\
z_2 \\
z_3
\end{array}\right).\end{eqnarray}
From this we can easily determine the lines corresponding to $M.$ 
Indeed, the $3\times 4$ matrix in~(\ref{enacba7}) has rank 2 at points 
$$(1,0,0),\ (0,1,0),\ (0,0,1),\ (1,1,1),\ (1,\omega,\omega^2),\ (1,\omega^2,\omega)$$
in $\PP^2$, which induce mutually skew lines 
$${\omega\,1\,0\,0\choose 0\,0\,\omega\,1},\ {1\,1\,0\,0\choose 0\,0\,1\,\omega},\ {1\,\omega\,0\,0\choose 0\,0\,1\,1},\ 
{1\,0\,0\,1\choose 0\,1\,1\,0},\ {1\,0\,0\,\omega\choose 0\,1\,\omega\,0},\ {\omega\,0\,0\,1\choose 0\,\omega\,1\,0}.$$
On the other hand $M^t$ induces lines 
$${1\,1\,0\,0\choose 0\,0\,1\,1},\ {1\,\omega\,0\,0\choose 0\,0\,\omega\,1},\ {\omega\,1\,0\,0\choose 0\,0\,1\,\omega},\ 
{\omega\,0\,0\,1\choose 0\,1\,\omega\,0},\ {1\,0\,0\,1\choose 0\,\omega\,1\,0},\ {1\,0\,0\,\omega\choose 0\,1\,1\,0}$$
over points
$$(1,0,0),\ (0,1,0),\ (0,0,1),\ (\omega^2,1,1),\ (1,1,\omega^2),\ (1,\omega^2,1).$$
Note that the above sets of lines form a double-six.\\

Consider next
$$M'=\left(\begin{array}{ccc}
0 & z_2+\omega z_3 & z_0+\omega z_1 \\
\omega z_0+z_1 & 0 & z_2+z_3  \\
\omega z_2+z_3  & z_0+z_1  & 0
\end{array}\right)$$ 
The columns of $M'$ determine lines 
$${\omega\,1\,0\,0\choose 0\,0\,\omega\,1},{1\,1\,0\,0\choose 0\,0\,1\,\omega},{1\,\omega\,0\,0\choose 0\,0\,1\,1}$$
and the rows of $M'$ (which are the columns of $(M')^t$) determine lines
$${1\,\omega\,0\,0\choose 0\,0\,1\,\omega},{\omega\,1\,0\,0\choose 0\,0\,1\,1},{1\,1\,0\,0\choose 0\,0\,\omega\,1}$$
which together form a half of a double-six.
It follows that the double-sixes corresponding to $M$ and $M'$ are different and so, they are not equivalent.\\

In the same way we show that the half  of a double-six corresponding to
$$M''=\left(\begin{array}{ccc}
0 & \omega z_0+z_1+z_2+\omega z_3 & \omega z_0+ z_1+\omega^2 z_2+z_3 \\
\omega z_0+\omega^2 z_1+z_2+z_3 & 0 & \omega z_0+z_3  \\
z_0+\omega z_1+\omega z_2+\omega z_3  & z_1+z_2  & 0
\end{array}\right)$$
is equal to
$${1\,\omega\,0\,0\choose 0\,0\,1\,1},\ {1\,0\,0\,1\choose 0\,1\,1\,0},\ 
{\omega\,0\,0\,1\choose 0\,\omega\,1\,0}\ \mbox{ and }\ 
{\omega\,1\,0\,0\choose 0\,0\,1\,\omega},\ {\omega\,0\,0\,1\choose 0\,1\,\omega\,0},\ {1\,0\,0\,\omega\choose 0\,1\,1\,0}.
$$
This set uniquely extends to the double-six corresponding to $M.$ Therefore $M$ and $M''$ are equivalent.
}
\end{example}

\section{Self-adjoint determinantal representations}
\label{selfadjsec}

From now on let the ground field $k$ be the complex field $\CC$. 
In Corollary~\ref{symm} we saw that a smooth cubic surface cannot have symmetric determinantal representations.
It is natural to ask whether a real surface $S$ (i.e., a surface such that all the coefficients of the defining 
homogeneous polynomial $F$ are real) admits a self-adjoint determinantal representation.
We will consider a determinantal representations $U(z_0,z_1,z_2,z_3)$ satisfying
$$\mbox{det}U=c F\ \mbox{ for some }\ c\in\CC,\ c\neq 0$$
and
$$\begin{array}{ccl}
U(z_0,z_1,z_2,z_3) &=&\left(u_{ij}^0 z_0 +u_{ij}^1 z_1 +u_{ij}^2 z_2 +u_{ij}^3 z_3\right)_{i,j}\\
\| & & \\
U^{\ast}(z_0,z_1,z_2,z_3)&=&\left(\overline{u_{ji}^0} z_0 +\overline{u_{ji}^1} z_1 +\overline{u_{ji}^2} z_2 +
\overline{u_{ji}^3} z_3\right)_{i,j}.
\end{array}$$

Assume that $U=U^{\ast}$. By Proposition~\ref{mainth} there exists an equivalent representation 
$$\Re=X U Y,\ \ \ X,Y\in{\GL}_3(\CC)$$ with zeros along the diagonal. Note that $\Re$ and $\Re^{\ast}$ 
are equivalent representations since
$$\Re^{\ast}=Y^{\ast}U X^{\ast}=Y^{\ast}X^{-1}\Re Y^{-1}X^{\ast}$$
and thus they correspond to the same six skew lines. 
We denote by $a_1\ldots a_6 \choose b_1\ldots b_6$ the double-six corresponding to
$\Re$ and $\Re^t$. 
The exact relations are the following:
$$\begin{array}{c}\Re=\\
 \\
 \\
\end{array}
\begin{array}{c}
\left(\begin{array}{c}
0   \\
\pi_{21}\\
\pi_{31}
\end{array}\right. \\
\downarrow \\
a_1 
\end{array}
\begin{array}{c}
\begin{array}{c}
\pi_{12}   \\
0\\
\pi_{32}
\end{array} \\
\downarrow \\
a_2 
\end{array}
\begin{array}{c}
\left.\begin{array}{c}
\pi_{13}  \\
\pi_{23}\\
0
\end{array}\right) \\
\downarrow \\
a_3 
\end{array}
\begin{array}{c}
\begin{array}{c}
\rightarrow   \\
\rightarrow  \\
\rightarrow 
\end{array} \\
\\
\\
\end{array}
\begin{array}{c}
\begin{array}{c}
b_1   \\
b_2 \\
b_3
\end{array} \\
\\
\\
\end{array}
$$
and
$$\begin{array}{c}\Re^{\ast}=\\
 \\
 \\
 \\
 \\
\end{array}
\begin{array}{c}
\left(\begin{array}{c}
0   \\
\overline{\pi_{12}}\\
\overline{\pi_{13}}
\end{array}\right. \\
\downarrow \\
a_i \\
\| \\
\overline{b_1}
\end{array}
\begin{array}{c}
\begin{array}{c}
\overline{\pi_{21}}   \\
0\\
\overline{\pi_{23}}
\end{array} \\
\downarrow \\
a_j \\
\| \\
\overline{b_2}
\end{array}
\begin{array}{c}
\left.\begin{array}{c}
\overline{\pi_{31}}  \\
\overline{\pi_{32}}\\
0
\end{array}\right) \\
\downarrow \\
a_k\\
\| \\
\overline{b_3} 
\end{array}
\begin{array}{c}
\begin{array}{c}
\rightarrow   \\
\rightarrow  \\
\rightarrow 
\end{array} \\
\\
\\
\\
\\
\end{array}
\begin{array}{c}
\begin{array}{c}
b_i=\overline{a_1}    \\
b_j=\overline{a_2} \\
b_k=\overline{a_3}
\end{array} \\
\\
\\
\\
\\
\end{array}
$$
for some $i\neq j\neq k\in\{1,\ldots,6\}$.
This gives a necessary condition for $S$ to have a self-adjoint determinantal representation. In fact, we
will prove that it is sufficient as well.

\begin{lemma} \label{selfadd}
A real smooth cubic surface $S$ has a self-adjoint determinantal representation 
if and only if there exists a double-six 
$$\left(\begin{array}{ccc}
a_1 & \ldots & a_6 \\
b_1 & \ldots & b_6
\end{array}\right)$$
and $i\neq j\neq k\in\{1,\ldots,6\}$
such that
$$\begin{array}{ccc}
a_i=\overline{b_1}, &  a_j=\overline{b_2}, & a_k=\overline{b_3}, \\
b_i=\overline{a_1}, &  b_j=\overline{a_2}, & b_k=\overline{a_3}
\end{array}.$$
\end{lemma}

\begin{proof}
It remains to show that the condition is sufficient. From the given double-six we
construct a determinantal representation $\Re$ with zeros along the
diagonal, such that its columns induce $a_1,a_2,a_3$ and its rows induce $b_1,b_2,b_3$. 
Then our assumption implies that $\Re^{\ast}$ induces the same six skew lines 
and is thus equivalent to $\Re$.
This means 
\begin{equation}\label{itisdiag} \Re^{\ast}S=T \Re \end{equation}
for some $S,\ T \in\GL_3(\CC).$ 

Recall that the intersections of planes given by the columns of $\Re$  
$$\Re   \left(\begin{array}{c}
1\\
0 \\
0
\end{array}\right),\ \Re \left(\begin{array}{c}
0\\
1 \\
0
\end{array}\right),\ \Re \left(\begin{array}{c}
0\\
0 \\
1
\end{array}\right)$$
induce the lines $a_1,a_2,a_3$, respectively. Since $\Re^{\ast}$ is equivalent to $\Re$ there exist uniquely determined
points 
$$(\zeta_i, \eta_i, \xi_i)\in \PP^2, \ \ i=1,2,3,$$
such that the intersections of planes corresponding to
$$\Re^{\ast}  \left(\begin{array}{c}
\zeta_1\\
\eta_1 \\
\xi_1
\end{array}\right),\ \Re^{\ast}  \left(\begin{array}{c}
\zeta_2\\
\eta_2 \\
\xi_2
\end{array}\right),\ \Re^{\ast}  \left(\begin{array}{c}
\zeta_3\\
\eta_3 \\
\xi_3
\end{array}\right)$$
also induce $a_1,a_2,a_3$, respectively. Evaluate (\ref{itisdiag}) on the vectors of the standard basis to prove that 
$$S= \left(\begin{array}{ccc}
\alpha_1 \zeta_1 & \alpha_2 \zeta_2 & \alpha_3 \zeta_3 \\
\alpha_1 \eta_1 & \alpha_2 \eta_2 & \alpha_3 \eta_3 \\
\alpha_1 \xi_1 & \alpha_2 \xi_2 & \alpha_3 \xi_3
\end{array}\right)$$
for some nonzero $\alpha_1,\alpha_2,\alpha_3\in\CC$.

Next observe that $\overline{\Re}$ over the standard vectors
induces the lines $\overline{a_1},\overline{a_2},\overline{a_3}$, respectively. In the same way
$\overline{\Re^{\ast}}=\Re^t$ over 
$(\overline{\zeta_i}, \overline{\eta_i}, \overline{\xi_i}),\ i=1,2,3$ again induces 
$\overline{a_1},\overline{a_2},\overline{a_3}$.
Consider the transpose of (\ref{itisdiag})
$$\Re^t T^t=S^t (\Re^{\ast})^t=S^t \overline{\Re}$$
This equality evaluated on the vectors of the standard basis implies that
$$T^t= \left(\begin{array}{ccc}
\beta_1 \overline{\zeta_1} & \beta_2 \overline{\zeta_2} & \beta_3 \overline{\zeta_3} \\
\beta_1 \overline{\eta_1} & \beta_2 \overline{\eta_2} & \beta_3 \overline{\eta_3} \\
\beta_1 \overline{\xi_1} & \beta_2 \overline{\xi_2} & \beta_3 \overline{\xi_3}
\end{array}\right)$$
for some nonzero $\beta_1,\beta_2,\beta_3\in\CC$.

Thus we proved that 
$$S=T^{\ast} D,\ \mbox{ where }\ D=
\left(\begin{array}{ccc}
\frac{\alpha_1}{\overline{\beta_1}} & 0 & 0 \\
0 & \frac{\alpha_2}{\overline{\beta_2}} & 0 \\
0 & 0 & \frac{\alpha_3}{\overline{\beta_3}}
\end{array}\right).$$
Then (\ref{itisdiag}) becomes
$$T \Re=\Re^{\ast}T^{\ast} D=(T \Re)^{\ast} D=\overline{D} T \Re D.$$
Since $T \Re$ is a determinantal representation of an irreducible surface, at most one entry in each row or column 
is zero. This implies that $D=\gamma\,$Id
for some $|\gamma|=1$
and the representation
$$\sqrt{\overline{\gamma}}\,T \Re =\sqrt{\overline{\gamma}}\, (T \Re)^{\ast} D
=\sqrt{\overline{\gamma}}\, \gamma\, (T \Re)^{\ast}=\sqrt{\gamma}\, (T \Re)^{\ast}$$
is self-adjoint.
\end{proof}

\begin{remark}\rm{ Note that a line $l$ lies on a real surface $S$ if and only if $\overline{l}$ lies on $S$.
There are three types of lines on $S$: 

\noindent 1){\it Conjugate of the 1st kind:}
\begin{itemize} 
\item $l\cap\overline{l}=$ one real point if and only if
\item $l$ has exactly 1 real point and $l\neq\overline{l}$.
\end{itemize}
2) {\it Conjugate of the 2nd kind:}
\begin{itemize} 
\item $l\cap\overline{l}=\emptyset$ if and only if
\item $l$ has no real points (i.e., no points which has all components real).
\end{itemize}
3) {\it Real line:}
\begin{itemize} 
\item $l=\overline{l}$ if and only if
\item $l$ is a real line.
\end{itemize}

Suppose that 
$\left(\begin{array}{ccc}
a_1 & \ldots & a_6 \\
b_1 & \ldots & b_6
\end{array}\right)$
is a double-six that satisfies the conditions of Lemma~\ref{selfadd}. Then the lines $a_1,a_2$ and $a_3$ 
cannot be real lines since all the $12$ lines of a double-six are distinct. If $a_1$ is a conjugate of the 
2nd kind then $i=1$ and 
$\overline{a_1}=b_1$. On the other hand, if $a_1$ is a conjugate of the 1st kind, then $i\neq 1$.
In this case $\pi_{i1}=\Span{a_1,\overline{a_1}}$ and $\pi_{1i}=\Span{a_i,\overline{a_i}}$ are real planes 
(i.e., given by a linear polynomial with real coefficients) and $c_{1i}=\pi_{i1}\cap\pi_{1i}$ is 
a real line.}
\end{remark}

A double-six $\left(\begin{array}{ccc}
a_1 & \ldots & a_6 \\
b_1 & \ldots & b_6
\end{array}\right)$
is called {\em mutually self-conjugate} if 
$\left\{b_1,\ldots,b_6\right\}=\left\{\overline{a_1},\ldots,\overline{a_6}\right\}$ as sets.
Using the fact that a half of the double-six uniquely determines the whole double-six we get the following:
\begin{corollary} 
A real smooth cubic surface $S$ has a self-adjoint determinantal representation 
if and only if there exists a mutually self-conjugate double-six on $S$.
\end{corollary}

Next we consider all nonequivalent self-adjoint representations of a given real cubic surface $S$. 
In Segre~\cite[pp. 40-55]{segre} we find the following classification:
\begin{theorem} A nonsingular real cubic surface can only be one of the 5 types specified in the table:\\

\begin{tabular}{||c|ccc||}\hline
Type                &  & Number of lines: &  \\
   & Real &  1st kind &  2nd kind \\ \hline
$F_1$ & 27 & 0 & 0 \\ 
$F_2$ & 15 & 0 & 12 \\
$F_3$ & 7 & 4 & 16 \\
$F_4$ & 3 & 12 & 12 \\
$F_5$ & 3 & 24 & 0 \\ 
\hline
\end{tabular}\\

After a suitable permutation of indexes, a mutually self-conjugate double-six is one of the
following $4$ types: a double-six {\it of the I-st kind } is of the form
$$\left(\begin{array}{cccccc}
a_1 & a_2&a_3&a_4&a_5 & a_6 \\
\overline{a_1} & \overline{a_2}&\overline{a_3}&\overline{a_4}&\overline{a_5} & \overline{a_6}
\end{array}\right),$$
a double-six {\it of the II-nd kind } is of the form
$$\left(\begin{array}{cccccc}
a_1 & a_2&a_3&a_4&a_5 & a_6 \\
\overline{a_2} & \overline{a_1}&\overline{a_3}&\overline{a_4}&\overline{a_5} & \overline{a_6}
\end{array}\right),$$
a double-six {\it of the III-rd kind } is of the form
$$\left(\begin{array}{cccccc}
a_1 & a_2&a_3&a_4&a_5 & a_6 \\
\overline{a_2} & \overline{a_1}&\overline{a_4}&\overline{a_3}&\overline{a_5} & \overline{a_6}
\end{array}\right)$$
and a double-six {\it of the IV-th kind } is of the form
$$\left(\begin{array}{cccccc}
a_1 & a_2&a_3&a_4&a_5 & a_6 \\
\overline{a_2} & \overline{a_1}&\overline{a_4}&\overline{a_3}&\overline{a_6} & \overline{a_5}
\end{array}\right).$$

Moreover, all mutually self-conjugate double-sixes are specified by\\

\begin{tabular}{||c|ccccc||}\hline
Type            & $F_1$ & $F_2$ & $F_3$ & $F_4$ & $F_5$\\  \hline
Number$_{Kind}$    &           $0$ & $1_{I}$& $2_{II}$ & $3_{III}$ & $12_{IV}$\\ 
             \hline
\end{tabular}
\end{theorem}

The four kinds of mutually self-conjugate double-sixes were introduced by Cremona \cite{cre1}. 
 
Every mutually self-conjugate double-six induces two nonequivalent self-adjoint determinantal representations. As 
a consequence we obtain the following result.
\begin{corollary} A real cubic surface has the following number of nonequivalent self-adjoint determinantal 
representations:\\

\begin{tabular}{||c|ccccc||}\hline
Type of the surface          & $F_1$ & $F_2$ & $F_3$ & $F_4$ & $F_5$\\  \hline
Number of s.a. reps    &           $0$ & $2$& $4$ & $6$ & $24$\\ 
             \hline
\end{tabular}
\end{corollary}

\begin{definition}{\rm Two self-adjoint determinantal representations $U_1,\ U_2$ of $S$ are 
{\it hermitean equivalent} if there exists a complex matrix $X\in \GL_3(\CC)$ such that
$$U_2=X U_1 X^{\ast}.$$}
\end{definition}

\begin{proposition}  Two self-adjoint determinantal representations $U_1,\ U_2$ of $S$ are 
equivalent if and only if $U_1$ is hermitean equivalent to either $U_2$ or $-U_2$.
\end{proposition}

\begin{proof} Assume that $U_1$ and $U_2$ are equivalent. Then there exist $X, Y\in \GL_3(\CC)$ such
that $U_2=X U_1 Y.$ 

Observe that $$X U_1 X^{\ast}=X U_1 Y Y^{-1} X^{\ast}=U_2 Y^{-1} X^{\ast}$$
is self-adjoint. Denote $Y^{-1} X^{\ast}$ by $Z$. Then 
\begin{equation}\label{zismultid} U_2 Z=Z^{\ast} U_2.\end{equation} 
Recall that there exist 6 unique 
points in general position $(\zeta_i,\eta_i,\xi_i)\in\PP^2,\ i=1,\ldots ,6$, such that for every $i$ equations 
$$U_2\cdot \left(\begin{array}{c}
\zeta_i\\
\eta_i \\
\xi_i
\end{array}\right)=0$$
define 
three planes which intersect in the line $a_i$. 
We can conclude from  (\ref{zismultid}) that $Z$ is a real multiple of 
the identity matrix. Indeed, for $i=1,\ldots ,6$ evaluate (\ref{zismultid})  on $(\zeta_i,\eta_i,\xi_i)$, 
which shows that 
$$Z \left(\begin{array}{c}
\zeta_i\\
\eta_i \\
\xi_i
\end{array}\right)=\lambda_i\left(\begin{array}{c}
\zeta_i\\
\eta_i \\
\xi_i
\end{array}\right)\ \mbox{ for some }\ 0\neq \lambda_i\in\CC.$$
Because the 6 points are in general position, every subset of 3 points generates a basis for $\CC^3$. 
Therefore all $\lambda_i$ are equal. Moreover, (\ref{zismultid}) implies that $\lambda_i$ is real and thus
$Z$ is a real multiple of Id.

This means that $Y=k X^{\ast}$ for some $k\in\RR$ and 
$$\begin{array}{ccc}
U_2=\sqrt{k}\, X\, U_1\, \sqrt{k}\,X^{\ast} & \mbox{ if } & k>0, \\
 & & \\
-U_2=\sqrt{-k}\, X\, U_1\, \sqrt{-k}\,X^{\ast} & \mbox{ if } & k<0,
\end{array}$$
which concludes the proof.
\end{proof}

\begin{example}\label{exexex}{\rm Lines corresponding to $M'$ in Example~\ref{exferm} satisfy the conditions of 
Lemma~\ref{selfadd}. Indeed $M'$ and $(M')^t$ induce the double-six
$$\left(\begin{array}{c}
{\omega\,1\,0\,0\choose 0\,0\,\omega\,1},{1\,1\,0\,0\choose 0\,0\,1\,\omega},{1\,\omega\,0\,0\choose 0\,0\,1\,1},
{1\,0\,0\,1\choose 0\,1\,\omega\,0},{\omega\,0\,0\,1\choose 0\,1\,1\,0},{1\,0\,0\,\omega\choose 0\,\omega\,1\,0}\\
\\
{1\,\omega\,0\,0\choose 0\,0\,1\,\omega},{\omega\,1\,0\,0\choose 0\,0\,1\,1},{1\,1\,0\,0\choose 0\,0\,\omega\,1},
{1\,0\,0\,\omega\choose 0\,1\,1\,0},{1\,0\,0\,1\choose 0\,\omega\,1\,0},{\omega\,0\,0\,1\choose 0\,1\,\omega\,0}
\end{array}\right)$$
which is of the $III$-rd kind.

Then 
$$U_1=\left(\begin{array}{ccc}
\omega^2 & 0 & 0 \\
0 & 0 & 1 \\
0 & 1 & 0
\end{array}\right)
\cdot M'=
\left(\begin{array}{ccc}
0 & \omega^2 z_2+z_3 & \omega^2 z_0+z_1 \\
\omega z_2+z_3 &  z_0+ z_1 & 0 \\
\omega z_0+z_1 & 0 & z_2+ z_3
\end{array}\right)$$ 
is self-adjoint.
Note that the representation
$$N=\left(\begin{array}{ccc}
0 & z_0+z_1+ z_2+\omega z_3 & \omega^2 z_0+z_1+z_2+z_3 \\
z_0+\omega^2 z_1+ z_2+z_3 & 0 & z_2+ z_3 \\
 z_0+z_1+\omega z_2+z_3 & \omega^2(z_0+z_1) & 0
\end{array}\right)$$
induces the same skew lines and is thus equivalent to $M'$. 
Also
$$\begin{array}{rl}U_2 & =\left(\begin{array}{ccc}
\omega & 2+\omega^2 & 1-\omega^2 \\
0 & 0 & \omega \\
0 & 1 & 0
\end{array}\right)
\cdot N \\
 & \\
 & =\left(\begin{array}{ccc}
3(z_0+z_3) & \omega^2 z_0+\omega^2 z_1+ \omega z_2+\omega^2 z_3 & z_0+\omega z_1+z_2+z_3 \\
\omega z_0+\omega z_1+ \omega^2 z_2+\omega z_3 &  z_0+ z_1 & 0\\
z_0+\omega^2 z_1+z_2+z_3 & 0 & z_2+ z_3
\end{array}\right)
\end{array}$$ 
is self-adjoint.
It is easy to check that
$$U_1=\left(\begin{array}{ccc}
1 & -\omega^2 & -1 \\
0 & \omega^2 & 0 \\
0 & 0 & \omega
\end{array}\right)\cdot U_2\cdot \left(\begin{array}{ccc}
1&0&0 \\
-\omega & \omega & 0 \\
-1 & 0 & \omega^2
\end{array}\right).$$

From the list of 27 lines on Fermat surface in Example \ref{exferm} or the fact that it has a conjugate
double six of the $III$-rd kind we conclude that it is of the Segre
type $F_4$ and thus it has $6$ nonequivalent self-adjoint determinantal representations. Observe that determinantal
representations $M$ of (\ref{Fermat1}) is not equivalent to a self-adjoint one.}
\end{example}

\section{Definite determinantal representations}

In this section we consider the following question: when do the $4$ coefficient matrices of a self-adjoint determinantal
representation of a cubic surface admit a real linear combination which is positive definite?

\begin{definition} \rm A self-adjoint determinantal representation $U=z_0 U_0+z_1 U_1+z_2 U_2+z_3 U_3$
of $S$ is called \textit{definite} if there exist $c_0,\,c_1,\,c_2,\,c_3\in \RR$ such that
$$c_0 U_0+c_1 U_1+c_2 U_2+c_3 U_3>0$$
and is \textit{indefinite} otherwise.

A nonzero vector $h\in \CC^3$ is \textit{self-orthogonal vector} of $U$ if
$$h^{\ast} U_i h=0,\ \mbox{ for all }i=0,1,2,3.$$
\end{definition}

Note that the notions defined above are projectively invariant. Furthermore, the definiteness is preserved under hermitean 
equivalence of determinantal representations. Clearly, a determinantal representation with a self-orthogonal vector is 
indefinite.

\begin{theorem} Every self-adjoint determinantal representation of a surface of types $F_2,F_3$ and $F_4$ 
has a self-orthogonal vector and is therefore indefinite.
\end{theorem}

\begin{proof} In Section~\ref{selfadjsec} we saw that self-adjoint determinantal representations are induced by mutually
self-conjugate double-sixes. Every double-six of a surface of type $F_2$ is of the $I-$st kind and the corresponding 
representation of type $\Re$ is of the form
$$\begin{array}{c}
\left(\begin{array}{c}
0   \\
\overline{\pi_{12}}\\
\overline{\pi_{13}}
\end{array}\right. \\
\downarrow \\
a_1 
\end{array}
\begin{array}{c}
\begin{array}{c}
\pi_{12}   \\
0\\
\overline{\pi_{23}}
\end{array} \\
\downarrow \\
a_2 
\end{array}
\begin{array}{c}
\left.\begin{array}{c}
\pi_{13}  \\
\pi_{23}\\
0
\end{array}\right) \\
\downarrow \\
a_3 
\end{array}
\begin{array}{c}
\begin{array}{c}
\rightarrow   \\
\rightarrow  \\
\rightarrow 
\end{array} \\
\\
\\
\end{array}
\begin{array}{c}
\begin{array}{c}
\overline{a_1}   \\
\overline{a_2} \\
\overline{a_3}
\end{array} \\
\\
\\
\end{array}
$$
Since, $\pi_{21}=\Span{a_1,\overline{a_2}}=\overline{\Span{a_2,\overline{a_1}}}=\overline{\pi_{12}}$, and similarly
$\pi_{31}=\overline{\pi_{13}}$ and $\pi_{32}=\overline{\pi_{23}}$, it follows that the representation is self-adjoint
and each of the coordinate vectors, e.g., $\left(\begin{array}{ccc}1&0&0\end{array}\right)^*$, is self-orthogonal.

Double-sixes of the $II-$nd and $III-$rd kind are, after a reordering of lines in the latter case, of the form
$$\left(\begin{array}{cccc}
a_1 & a_2&a_3&\cdots \\
\overline{a_2} & \overline{a_1}&\overline{a_3}&\cdots
\end{array}\right).$$
The corresponding representation of type $\Re$ and an equivalent representation are 
$$\begin{array}{c}
\left(\begin{array}{c}
0   \\
\pi_{21}\\
\pi_{31}
\end{array}\right. \\
\downarrow \\
a_1 
\end{array}
\begin{array}{c}
\begin{array}{c}
\pi_{12}   \\
0\\
\pi_{32}
\end{array} \\
\downarrow \\
a_2 
\end{array}
\begin{array}{c}
\left.\begin{array}{c}
\pi_{13}  \\
\pi_{23}\\
0
\end{array}\right) \\
\downarrow \\
a_3 
\end{array}
\begin{array}{c}
\begin{array}{c}
\rightarrow   \\
\rightarrow  \\
\rightarrow 
\end{array} \\
\\
\\
\end{array}
\begin{array}{c}
\begin{array}{c}
\overline{a_2}   \\
\overline{a_1} \\
\overline{a_3}
\end{array} \\
\\
\\
\end{array}
\begin{array}{c}\\
 \sim \\
 \\
 \\
 \\
\end{array}
\begin{array}{c}
\left(\begin{array}{c}
\pi_{21}   \\
0\\
\pi_{31}
\end{array}\right. \\
\downarrow \\
a_1 
\end{array}
\begin{array}{c}
\begin{array}{c}
0   \\
\pi_{12}\\
\pi_{32}
\end{array} \\
\downarrow \\
a_2 
\end{array}
\begin{array}{c}
\left.\begin{array}{c}
\pi_{23}  \\
\pi_{13}\\
0
\end{array}\right) \\
\downarrow \\
a_3 
\end{array}
\begin{array}{c}
\begin{array}{c}
\rightarrow   \\
\rightarrow  \\
\rightarrow 
\end{array} \\
\\
\\
\end{array}
\begin{array}{c}
\begin{array}{c}
\overline{a_1}   \\
\overline{a_2} \\
\overline{a_3}
\end{array} \\
\\
\\
\end{array}
$$
The latter representation is self-adjoint since $\pi_{21}=\Span{a_1,\overline{a_1}}$ and 
$\pi_{12}=\Span{a_2,\overline{a_2}}$ are real 
planes and
$$\pi_{31}=\Span{a_1,\overline{a_3}}=\overline{\Span{a_3,\overline{a_1}}}=\overline{\pi_{23}},\ \ \
\pi_{32}=\Span{a_2,\overline{a_3}}=\overline{\Span{a_3,\overline{a_2}}}=\overline{\pi_{13}}.$$
Moreover, the vector $\left(\begin{array}{ccc}0&0&1\end{array}\right)^*$ is its self-orthogonal vector. 
Hence all the self-adjoint determinantal representations of a surface of types $F_2,F_3$ and $F_4$ are
indefinite.
\end{proof}

\begin{example} \rm Recall the self-adjoint representation of Fermat surface from Example~\ref{exexex}:
$$U=z_0\left(\begin{array}{ccc}
0&0&\omega^2 \\
0& 1 & 0 \\
\omega & 0 & 0
\end{array}\right)+
z_1\left(\begin{array}{ccc}
0&0&1 \\
0& 1 & 0 \\
1 & 0 & 0
\end{array}\right)+
z_2\left(\begin{array}{ccc}
0&\omega^2 &0 \\
\omega& 0 & 0 \\
0 & 0 & 1
\end{array}\right)+
z_3\left(\begin{array}{ccc}
0&1&0 \\
1& 0 & 0 \\
0 & 0 & 1
\end{array}\right)$$
It is easy to check that 
$\left(\begin{array}{ccc}\omega&0&0\end{array}\right)$
is a self-orthogonal vector of $U$.
\end{example}

It is natural to ask whether the coefficients of every indefinite determinantal representation of $S$ 
have a self-orthogonal vector. We will see in Example~\ref{exF5} that the answer is no.

\begin{theorem}\label{thmF5} Every real cubic surface of type $F_5$ has up to equivalence 16 definite determinantal representations 
(among the 24 nonequivalent self-adjoint representations).
\end{theorem}

Before we prove the theorem we recall some facts about the geometry of real cubic surfaces of type $F_5$
which are found in Segre~\cite{segre}.

The real part of a real cubic surface of type $F_5$ consists of two disconnected components, one of which is ovoidal.
On the other hand, the real part of real cubic surfaces of types $F_1,\ F_2,\ F_3$ and $F_4$ consist of a single connected component. 

We denote a real cubic surface of type $F_5$ by $F_5$ and we write $F_5(\RR)$ for its real part.
Then $F_5$ has $3$  real lines all of which are coplanar and lie in the non-ovoidal piece. 
Through every real line there are 4 additional tritangent planes, each 
intersecting $F_5$ in 2 conjugate lines of the 1st kind. Thus $F_5$ has $13=1+3\times 4$ real tritangent planes. 
The $24$ conjugate lines form $12$ self-conjugate double-sixes of the $IV-$th kind. 

A \textit{parabolic plane} through a real line $r$ is by definition the plane on which $r$ is tangent to the 
residual conic of intersection with the surface. The touching point is called a \textit{parabolic point} of $r$. 
Segre proves that every real line of a real cubic surface $F_5$ contains 2 real parabolic points.

Note that a tritangent plane intersects the real part of $F_5(\RR)$ in either three real lines or in a union of a real line and a point 
which is the real intersection point of the two conjugate complex lines of the tritangent plane. 
This point may lie on the real line. In fact, it lies on the real line if and only if it is an Eckardt point. (If 
all three lines on a tritagent plane meet in one point then such a point is called an {\it Eckardt point}.) 
  
Now fix a real line $r$ on $F_5$ and consider real planes through $r$.
Among those there are:
\begin{itemize}
\item 2 parabolic planes $\lambda,\ \lambda'$,
\item 2 tritangent planes $\beta,\ \beta'$ that intersect the ovoidal piece of $F_5(\RR)$,  
\item 2 tritangent planes $\gamma,\ \gamma'$ that do not intersect the ovoidal piece of $F_5(\RR)$.
\end{itemize}

\begin{figure}[ht] 
\begin{center}
\includegraphics*[5mm, 2cm][11cm, 9cm]{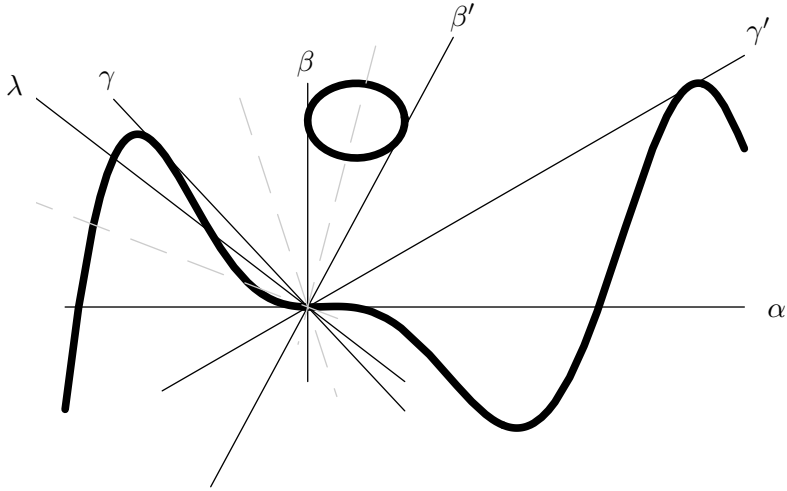}
  \put(-275,155){$\gamma$}
  \put(-310,151){$\lambda$}
  \put(-200,160){$\beta$}
  \put(-142,177){$\beta'$}
  \put(-30,172){$\gamma'$}
  \put(-22,67){$\alpha$}
\caption{$F_5$ surface}
\label{pict5}
\end{center}
\end{figure}

\begin{figure}[ht] 
\begin{center}
\includegraphics*[5mm, 4cm][11cm, 7cm]{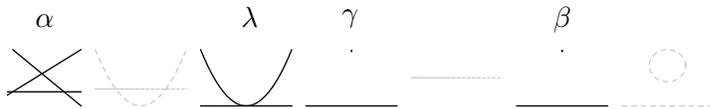}
  \put(-288,50){$\alpha$}
  \put(-210,50){$\lambda$}
  \put(-172,52){$\gamma$}
  \put(-92,50){$\beta$}
\caption{Planes containing $r$}
\label{pict6}
\end{center}
\end{figure}

Now we view our surface affinely and we consider an orientation on the pencil of real planes with axis $r$. 
We start with the plane $\alpha$ through the 3 real lines and fix a direction of rotation around $r$.
All the above planes are real and occur in the following order (see Figures \ref{pict5} and \ref{pict6}):
$$\alpha,\ \lambda,\ \gamma,\ \beta,\ \beta',\ \gamma',\ \lambda',\ \alpha.$$
Note that $\lambda$ and $\gamma$ may coincide. This happens if and only if the conjugate lines on $\gamma$ and $r$ intersect 
in an Eckardt point. The same holds for $\lambda'$ and $\gamma'$. No other planes can coincide. 

The ovoidal part of $F_5$ is wedged between $ \beta$ and $\beta'$ and the non-ovoidal between $\gamma'$ and $\gamma$ 
with respect to the above fixed orientation. On the other hand the interior of wedges between 
$\gamma,\ \beta$ and 
$\beta',\ \gamma'$ do not intersect $F_5(\RR)$.

\begin{proof}[Proof of Theorem~\ref{thmF5}] Fix a double six of the $IV$-th kind:
$$\left(\begin{array}{cccccc}
a_1 & a_2&a_3&a_4&a_5 & a_6 \\
\overline{a_2} & \overline{a_1}&\overline{a_4}&\overline{a_3}&\overline{a_6} & \overline{a_5}
\end{array}\right).$$
The corresponding determinantal representation of type $\Re$ is 
$$\begin{array}{c}\Re=\\
 \\
 \\
\end{array}
\begin{array}{c}
\left(\begin{array}{c}
0   \\
\pi_{21}\\
\pi_{31}
\end{array}\right. \\
\downarrow \\
a_1 
\end{array}
\begin{array}{c}
\begin{array}{c}
\pi_{12}   \\
0\\
\pi_{32}
\end{array} \\
\downarrow \\
a_2 
\end{array}
\begin{array}{c}
\left.\begin{array}{c}
\pi_{13}  \\
\pi_{23}\\
0
\end{array}\right) \\
\downarrow \\
a_3 
\end{array}
\begin{array}{c}
\begin{array}{c}
\rightarrow   \\
\rightarrow  \\
\rightarrow 
\end{array} \\
\\
\\
\end{array}
\begin{array}{c}
\begin{array}{c}
\overline{a_2}   \\
\overline{a_1} \\
\overline{a_4}
\end{array} \\
\\
\\
\end{array}.
$$
First we find a self-adjoint determinantal representation equivalent to $\Re$.
Planes $\pi_{12}$ and $\pi_{21}$ are real and intersect in a real line on $F_5$.
By Proposition~\ref{doublrem} there exists a point $(\eta,\xi,1)\in\PP^2$ such that 
the three planes corresponding to
$$\Re^t \left(\begin{array}{c}
\eta\\
\xi \\
1
\end{array}\right)$$
intersect in the line $\overline{a_3}$.
Then $U=$
$$
\begin{array}{c}\left(\begin{array}{ccc}
0&1&0\\
1&0&0\\
\eta & \xi&1
\end{array}\right)\Re=\\
 \\
 \\
\end{array}
\begin{array}{c}
\left(\begin{array}{c}
\pi_{21}   \\
0\\
\xi \pi_{21}+\pi_{31}
\end{array}\right. \\
\downarrow \\
a_1 
\end{array}
\begin{array}{c}
\begin{array}{c}
0   \\
\pi_{12}\\
\eta \pi_{12}+\pi_{32}
\end{array} \\
\downarrow \\
a_2 
\end{array}
\begin{array}{c}
\left.\begin{array}{c}
\pi_{23}  \\
\pi_{13}\\
\eta \pi_{13}+\xi \pi_{23}
\end{array}\right) \\
\downarrow \\
a_3 
\end{array}
\begin{array}{c}
\begin{array}{c}
\rightarrow   \\
\rightarrow  \\
\rightarrow 
\end{array} \\
\\
\\
\end{array}
\begin{array}{c}
\begin{array}{c}
\overline{a_1}   \\
\overline{a_2} \\
\overline{a_3}
\end{array} \\
\\
\\
\end{array}
$$
is self-adjoint since
$$\eta \pi_{12}+\pi_{32}=\Span{a_2,\overline{a_3}}=\overline{\pi_{13}},\ \
\xi \pi_{21}+\pi_{31}=\Span{a_1,\overline{a_3}}=\overline{\pi_{23}}$$
(after multiplying $\pi_{23}$ and $\pi_{13}$ by real constants if necessary) and 
$$\eta \pi_{13}+\xi \pi_{23}=\Span{a_3,\overline{a_3}}=\pi_{43}$$
is real.\\

To determine when the above determinantal representation $U$ is definite we introduce a system of 
quadratic forms:
$$\begin{array}{c}
\mathcal{Q}=(\overline{\alpha },\overline{\beta },\overline{\gamma })
\left(\begin{array}{ccc}
\pi_{21}&0&\pi_{23}\\
0&\pi_{12}&\pi_{13}\\
\overline{\pi_{23}}&\overline{\pi_{13}}&\pi_{43}
\end{array}\right)
\left(\begin{array}{c}
\alpha \\
\beta \\
\gamma 
\end{array}\right)=\\
\\
\frac{1}{\pi_{21}}(\alpha  \pi_{21}+\gamma  \pi_{23})(\overline{\alpha } \pi_{21}+\overline{\gamma } \overline{\pi_{23}})+
\frac{1}{\pi_{12}}(\beta  \pi_{12}+\gamma  \pi_{13})(\overline{\beta } \pi_{12}+\overline{\gamma } \overline{\pi_{13}})+\\
\\
\frac{1}{\pi_{12} \pi_{21}}(\pi_{12} \pi_{21} \pi_{43}- \pi_{12}\pi_{23}\overline{\pi_{23}} 
-\pi_{21} \pi_{13}\overline{\pi_{13}})\gamma  \overline{\gamma }.
\end{array}
$$
Observe that 
\begin{equation}\label{pipi}
\pi_{12} \pi_{21} \pi_{43}- \pi_{12}\pi_{23}\overline{\pi_{23}} 
-\pi_{21} \pi_{13}\overline{\pi_{13}}
\end{equation}
is the determinant of $U$ and therefore a defining equation of $F_5$.

Note that $U$ is definite if and only if $\mathcal{Q}$ contains a definite quadratic form, or equivalently, if there 
is a real point $P\in \PP^3$ such that $\pi_{12}(P),\pi_{21}(P)$
and $F_5(P)$ are all positive or all negative. 

Planes $\pi_{12}$ and $\pi_{21}$ divide $\PP^3$ into two wedges where $\pi_{12},\,\pi_{21}$ either have 
the same or different signs. We will prove that the non-ovoidal component of 
$F_5(\RR)$ lies in the wedge in which the values $\pi_{12}$ and $\pi_{21}$ have different signs. 
Consider the real line in $\pi_{43}\cap F_5(\RR)$  and choose a real point $P$ on it, that is not on $\pi_{13}$ or $\pi_{23}$.
Evaluate (\ref{pipi}) in $P$ to see that the values $\pi_{12}(P)$ and $\pi_{21}(P)$ must have different signs. 

Recall the geometry of real planes around the line $r$. It follows that the wedge where $\pi_{12},\,\pi_{21}$ have the same 
sign either contains the ovoidal piece, or the interior of the wedge does not intersect the real part of $F_5$. 
In the first case there exists a point in which $\pi_{12},\, \pi_{21}$ and (\ref{pipi}) have the same sign, which makes the
form $\mathcal{Q}(P)$ definite.

In the second case all the quadratic forms of $\mathcal{Q}$ are non-definite. Indeed, the plane $\pi_{43}$ 
cuts the wedge where $\pi_{12},\, \pi_{21}$ have the same sign. In the intersecting points (\ref{pipi}) has exactly 
the opposite sign to $\pi_{12},\, \pi_{21}$. Since the intersection with $F_5(\RR)$ is empty, the sign differs 
along the whole wedge.

After fixing a direction of rotation around $r$, the 4 tritangent planes through $r$ cut out 6 different wedges 
(i.e., different choices of $\pi_{12},\ \pi_{21}$)
in which the wedging two planes are either both positive or negative. From Figure \ref{pict5} we see that
4 wedges contain the ovoidal piece of $F_5(\RR)$ and thus induce positive-definite representations. 
The remaining 2 wedges induce non-definite representations.

Recall that every surface of type $F_5$ contains 12 mutually self-conjugate double-sixes of the $IV$-kind which induce
24 self-adjoint determinantal representations.
For a given
$$\left(\begin{array}{cccccc}
a_1 & a_2&a_3&a_4&a_5 & a_6 \\
\overline{a_2} & \overline{a_1}&\overline{a_4}&\overline{a_3}&\overline{a_6} & \overline{a_5}
\end{array}\right)$$ 
the wedges where the planes $\pi_{i,i+1}$ and $\pi_{i+1,i}$ $i=1,3,5,$  have the same sign
either all contain the ovoidal piece or their interior does not intersect $F_5(\RR)$. 
We have proved that $\frac{2}{3}$ of the wedges where both planes have the same sign contain the ovoidal piece. 
Thus it follows that out of $24$ self-adjoint determinantal 
representations $16$ are definite.
\end{proof}

\begin{example}\label{exF5} {\rm We use notation consistent with Proof of Theorem~\ref{thmF5}.
Let $F_5$ be a surface with equation 
$$\left(\frac{100}{24} z_0^2 + z_1^2\right)(z_0 + z_2) 
- z_3\left(z_3 - \frac{1}{2} z_2\right)\left(z_3 - \frac{2}{3} z_2\right)=0.$$
It has 3 real lines on the plane $\alpha:\ z_0 + z_2=0.$ Through the line 
$$r:\ z_0 + z_2=3 z_3-2 z_2=0$$
there are 4 real tritangent planes, each containing two intersecting complex conjugate lines:
$$\begin{array}{cl}
\gamma :&  z_0 + 0.98987 z_2 + 0.01519 z_3 , \\
\beta : &  z_0 + 0.01345 z_2 + 1.47982 z_3, \\
\beta': &  z_0 - 3.00333 z_2 + 6.00499 z_3, \\
\gamma': &  3 z_3-2 z_2  .
\end{array}$$
Determinantal representation 
$$\left(\begin{array}{ccc}
-z_0 - 0.98987 z_2 - 0.01519 z_3 & 0 & 
      \frac{ 2.04124 z_0 -  i z_1 + 8.14425 z_3}{28.68441(1 - i)}\\
  0& 3 z_3-2 z_2  & (1+i) (\frac{1}{2} z_0 -  i\frac{\sqrt{6}}{10} z_1)  \\
  \frac{ 2.04124 z_0 +i z_1 + 8.14425 z_3}{28.68441(1 +i)}&
(1-i) (\frac{1}{2} z_0 +  i\frac{\sqrt{6}}{10} z_1)  & -0.02020 z_3
\end{array}\right)$$
is definite since it is spanned by the wedge $\pi_{21}=\gamma,\ \pi_{12}=\gamma '$ which 
contains the ovoidal piece of $F_5$. For example, evaluate the representation at
$$z_0= 0.02,\ z_1= 0,\ z_2 = -1.2,\  z_3 = -0.3$$ and see that its eigenvalues
$1.50013,\ 1.17540,\ 0.00293$ are all positive.\\

On the other hand, determinantal representation 
$$U=\left(\begin{array}{ccc}
\frac{z_0 - 3.00333 z_2 + 6.00499 z_3}{-0.12000} & 0 & \frac{2.04124 z_0 - i z_1 - 0.00679 z_3}{0.49979(1-i)} \\
0 &  3 z_3-2 z_2  & (1+i) (\frac{1}{2} z_0 -  i\frac{\sqrt{6}}{10} z_1) \\
\frac{2.04124 z_0 + i z_1 - 0.00679 z_3}{0.49979(1+i)}  &  
(1-i) (\frac{1}{2} z_0 + i \frac{\sqrt{6}}{10} z_1) & 0.00666 z_3
\end{array}\right)$$
spanned by $\pi_{21}=\beta',\ \pi_{12}=\gamma '$ is non-definite.
In order to prove this, consider equivalent representation $U'=TUT^{\ast}$ where 
$$T=\left(\begin{array}{ccc}
0.02048(1+i) & 0.08309(1+i) & - i 0.24990 \\
-0.00666(1+i)&0.02040(1-i) & 12.25254\\
-0.24495(1+i)&1.00042(1+i) &0
\end{array}\right).$$
It is easy to check that 
$$\begin{array}{c}
U'=
z_3 \left(\begin{array}{ccc}
0&0&1\\
0&1&0\\
1&0&0
\end{array}\right)+
z_1 \left(\begin{array}{ccc}
0&1&0\\
1&0&0\\
0&0&0
\end{array}\right)+ \\
\\
z_2 \left(\begin{array}{ccc}
-0.00662 & -0.01360 i  & -0.58361 \\
 0.01361 i & 0.00055 & 0 \\
 -0.58361 & 0 & -1
\end{array}\right)+ 
z_0 \left(\begin{array}{ccc}
-0.09032 & 2.04691 i & 0.08361 \\
-2.04691 i & -0.16722 & 0.02718 i\\
0.08361 & -0.02719 i  & -1
\end{array}\right)
\end{array}$$
and Det$\,U=$Det$\,U'$.

With Mathematica 5.0 we compute the eigenvalues of $U'$. They are solutions 
$\lambda_i(z_0,z_1,z_2,z_3),\ i=1,2,3$ of the cubic equation
$$ \begin{array}{l}
\Lambda^3+\Lambda^2(1.25754 z_0 +  1.00607 z_2 - z_3)+\\
\Lambda (-3.92492 z_0^2 - z_1^2 + 0.41797 z_0 z_2 - 0.33472 z_2^2 - 1.25754 z_0 z_3 + 
0.16060 z_2 z_3 - z_3^2)  \\
-4.16667 z_0^3 - z_0 z_1^2 - 4.16667 z_0^2 z_2 - z_1^2 z_2 + 
      0.33333 z_2^2 z_3 - 1.16667 z_2 z_3^2 + z_3^3=0.
              \end{array}$$
From 
$$(\Lambda-\lambda_1)(\Lambda-\lambda_2)(\Lambda-\lambda_3)=\Lambda^3-\Lambda^2(\lambda_1+\lambda_2+\lambda_3)+
\Lambda(\lambda_1 \lambda_2+\lambda_1 \lambda_3+\lambda_2 \lambda_3)-\lambda_1 \lambda_2 \lambda_3$$
we observe the following: suppose that there exist values $z_0,z_1,z_2,z_3$ for which all $\lambda_i$ are of the same sign.
Then 
$$\begin{array}{l}
\lambda_1 \lambda_2+\lambda_1 \lambda_3+\lambda_2 \lambda_3=\\
-3.92492 z_0^2 - z_1^2 + 0.41797 z_0 z_2 - 0.33472 z_2^2 - 1.25754 z_0 z_3 + 
0.16060 z_2 z_3 - z_3^2
\end{array}$$ is positive. But this equals
$$-(1.87656 z_0 - 0.11137 z_2)^2 - (0.63516 z_0 + 0.98995 z_3)^2 -(0.56773 z_2 - 0.14144 z_3)^2-z_1^2$$
which leads to contradiction.
This finishes the proof that $U'$ is non-definite representation.\\

Finally we prove that (unlike surfaces of types $F_2,\ F_3,\ F_4$ ) the self-adjoint representation 
$U'=z_0 U_0+z_1 U_1+z_2 U_2+z_3 U_3$ 
has no self-orthogonal vector $h=(h_0,h_1,h_2)$. Again argue by contradiction. We can assume that $h_0=1$. 
From 
$$\begin{array}{l}
h U_3 h^{\ast}=h_1 \overline{h_1}+h_2+ \overline{h_2}=0 \\
h U_1 h^{\ast}=h_1 +\overline{h_1}=0
\end{array}$$
we get Re$(h_1)=0$ and Re$(h_2)=-\frac{1}{2}|h_1|^2$. Thus $h=(1,i k_1,-\frac{1}{2} k_1^2+i k_2)$ for some
$k_1,k_2\in\RR$.
Then
$$\begin{array}{l}
h U_2 h^{\ast}=-0.00662- 0.02722 k_1 + 0.58416 k_1^2 - 0.25 k_1^4 - k_2^2=0\\
h U_0 h^{\ast}=-0.09032+  4.09382 k_1- 0.25083 k_1^2 + 
0.02719 k_1^3 - 0.25 k_1^4 -k_2^2=0
\end{array}$$
It is easy to check that this system has no real solution. Therefore there are no self-orthogonal vectors 
for $U'$.}
\end{example}

\end{document}